# Solution of a NxN System of Linear algebraic Equations:

# 1 - The Steepest Descent Method Revisited

Dr. V. S. Patwardhan[1]

May 31, 2022

**Abstract**

This is the first in a series of papers which deal with the development of novel methods for solving a system of linear algebraic equations with a time complexity lower than existing algorithms. The NxN system of linear equations, **Ax** = **b**, is often solved iteratively by minimising the corresponding quadratic form using well known optimisation techniques. The simplest of these is the steepest descent method, whose approach to the solution is usually quite rapid at the beginning but slows down drastically after a few iterations. This paper investigates possible approaches which can reduce or avoid this slowing down. The two approaches used here involve random movement of the point between iterations, and possible matrix transformations between iterations. This paper reports the results of computational experiments and shows the remarkable improvement in performance of the steepest descent method that is possible. The approaches described here do not give a practical algorithm right away. However, they set the stage for developing practical algorithms based on SD alone, which will be presented in later publications.

## 1. Introduction

The solution of linear algebraic equations is a classical problem which is theoretically important and has many practical applications in engineering and science as well. The NxN system of linear equations can be written in a matrix form as **Ax** = **b**, where **A** is a NxN matrix, **b** is an N-vector, and **x** is the solution vector to be determined. The classical method involving gaussian elimination or other direct methods give exact solutions but are suitable only for moderately large problems since they have a time complexity of $O(N^3)$. In practical applications it is often unnecessary to get an exact solution. It is sufficient to get an approximate solution. Iterative methods essentially aim to start with a guess solution and improve it iteratively to get closer to the solution within acceptable accuracy. Methods of Jacobi, Gauss-Seidel, and Relaxation are well known but converge to the solution only if certain conditions are met. A quick summary of these methods can be found in standard books on linear algebra [for example, J. E. Gentle, 2007].

There are other methods such as the steepest descent method and the conjugate gradient method which are applicable when **A** is a symmetric positive definite matrix. This is not a great limitation, because the system of equations stated above, i.e. **Ax** = **b**, can be put in the form **A**$^T$**Ax** = **A**$^T$**b**, where the matrix **A**$^T$**A** is symmetric positive definite for any **A**. The steepest descent method and the conjugate gradient method are based on minimizing an appropriately defined quadratic function, using optimization techniques. Details of these methods, including convergence analysis, are available in standard books and reports [for example, J. R. Shewchuk, 1994]. The steepest descent method follows the path of steepest descent till a minimum point is reached along that line and repeats the procedure till the required convergence is achieved. The steepest descent method can be accelerated using the momentum concept to achieve faster convergence [Y. Nesterov, 1983; I. Sutskever et. Al., 2013]. The conjugate gradient method uses a different step size and converges exactly to the solution in N steps. It is however, often used as an iterative method since close approach to the solution is usually found in a much smaller number of steps.

____________________________________________________________________________

[1] Independent researcher. Formerly, Scientist G, National Chemical Laboratory, Pune 411008, India.
  Email : vspatw@gmail.com , URL : https://www.vspatwardhan.com

In this paper, we take a closer look at the steepest descent method (called as SD here). It is well known that SD converges to the solution in one step if the starting point **x** is located such that the line x-s happens to be one of the eigenvectors of the **A** matrix (i.e., one of the axes of the quadratic form). On the other hand, if the starting point **x** is located at the worst possible location, convergence can be so slow as to be impractical. However, if the starting point x is selected at random, the convergence shows an interesting behaviour. The objective of this paper is to investigate this behaviour through numerical experiments and set the stage for developing efficient algorithms which would be based on SD alone.

## 2. Generation of test problems

To study the behaviour of SD in greater detail, several problems were randomly generated. Practical problems usually involve highly sparse **A** matrices. However, results obtained with sparse problems are likely to be influenced by the sparsity structure. To avoid this source of variation, we only use dense **A** matrices here. The following procedure was used for generating the problem data, i.e. **A** and **b**. A random vector was generated with N(0,1) elements (i.e. standard normal deviates). The solution point **s** was chosen 10 units away from the origin in the direction of this random vector. The coefficient matrix **A** was then generated with N(0,1) elements. The **b** vector was generated as **b** = **As**. To select the initial point, another random vector was generated with N(0,1) elements, and the initial point $\mathbf{x}_0$ was selected by travelling from **s** in this random direction, a distance between 1 and 10 units. SD was then used starting from $\mathbf{x}_0$ and the progress of successive iterates $\mathbf{x}_k$ was monitored for k=1 to 10. The system size was varied from 10 to 1000. After each SD step, the point obtained was characterised by calculating (1) $\sigma_{Res}$ (i.e. the root mean square of the residuals for the current point x, defined as **r** = **b** − **Ax** ), (2) $d_{soln}$ (i.e. the distance from the solution point **s**), and (3) $\sigma_{\Delta x}$ (i.e. the root mean square value of the deviations ($\mathbf{x}_j$ - $\mathbf{s}_j$), for j = 1 to N). Two of these quantities are of course related as $d_{soln} = \sigma_{\Delta x} N^{0.5}$. During this procedure, the known solution point was used only for characterising the successive iterates as stated above. The SD algorithm itself for n steps, is shown below for clarity:

______________________________________________

*SD, Algorithm 1: Standard procedure*

Initial values:         **A**, **b**, and $\mathbf{x}_0$ = the starting point
For k = 1 to n:         $\mathbf{r}_k$ = **b** − **Ax**$_{k-1}$

$$\alpha_k = r_k^T r_k / ( r_k^T A r_k )$$
$$\mathbf{x}_k = \mathbf{x}_{k-1} + \alpha_k \mathbf{r}_k$$

___

## 3. Results of numerical experiments

Figure 1 shows results obtained for a specific problem for N = 10. The y-axis shows the normalized rms of residuals, $\sigma_{Res,k}/\sigma_{Res,init}$, where k is the iteration number. It is seen that the first iteration itself reduces $\sigma_{Res,k}$ by more than 60%, while further iterations reduce it further by much smaller amounts. This is typical of the progress obtained with SD, which seems to get almost stuck at some point.

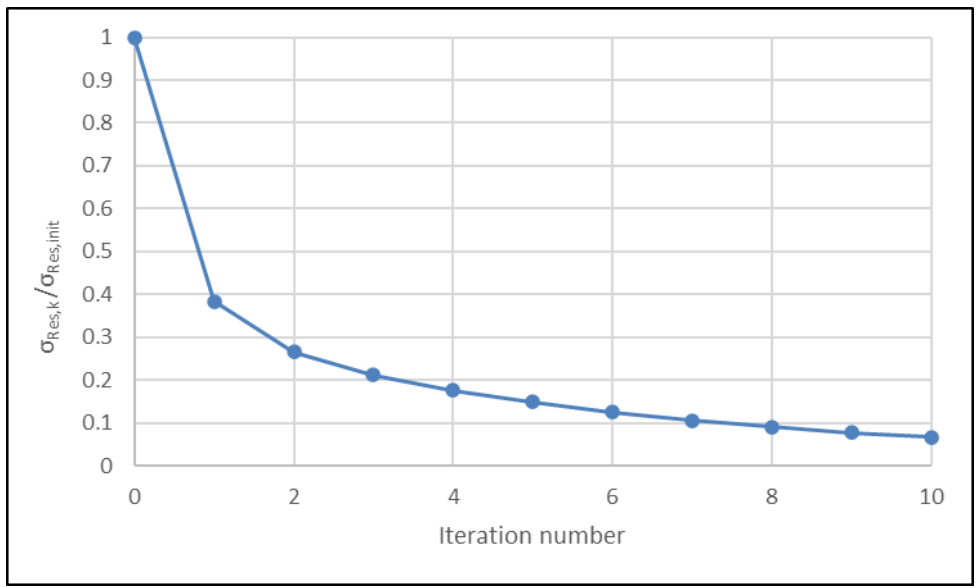

**Figure 1**. Residuals vs. iteration number for N = 10

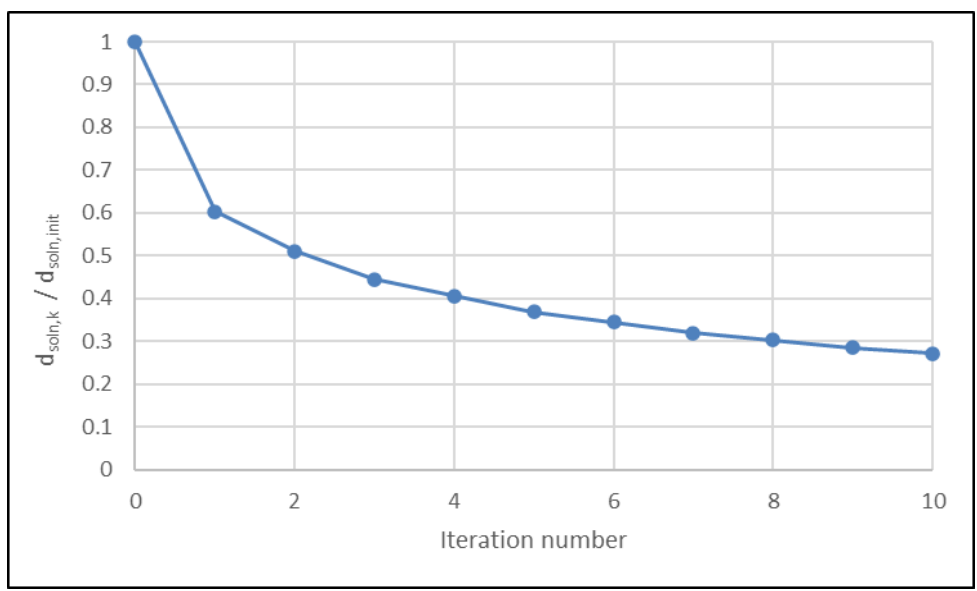

**Figure 2**. Distance from the solution vs. iteration number for N = 10

Figure 2 shows the normalized distance from the solution with iteration number. This figure shows that the distance from the solution reduces by 40% in one iteration itself, while further reduction slowly tapers off and further progress becomes slower and slower.

This kind of behaviour is typical of SD, where successive iterations give points which go in a zigzag manner, and the step size of each iteration goes on reducing rapidly. Figures 1 and 2 refer to one specific problem and a starting point. A different problem and a different starting point would give different numbers. To see the average behaviour of SD, many calculations were made with several starting points. The results for residuals, obtained for 10 different starting points are shown in Figure 3 for N = 10, and in Figure 4 for N = 1000.

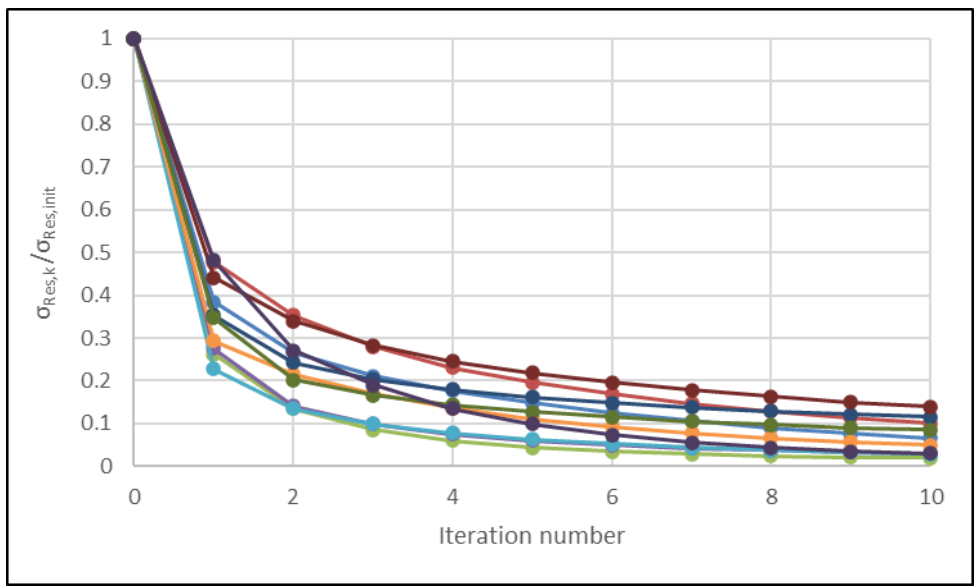

**Figure 3**. Residuals vs. iteration number for N = 10 for 10 values of $x_0$

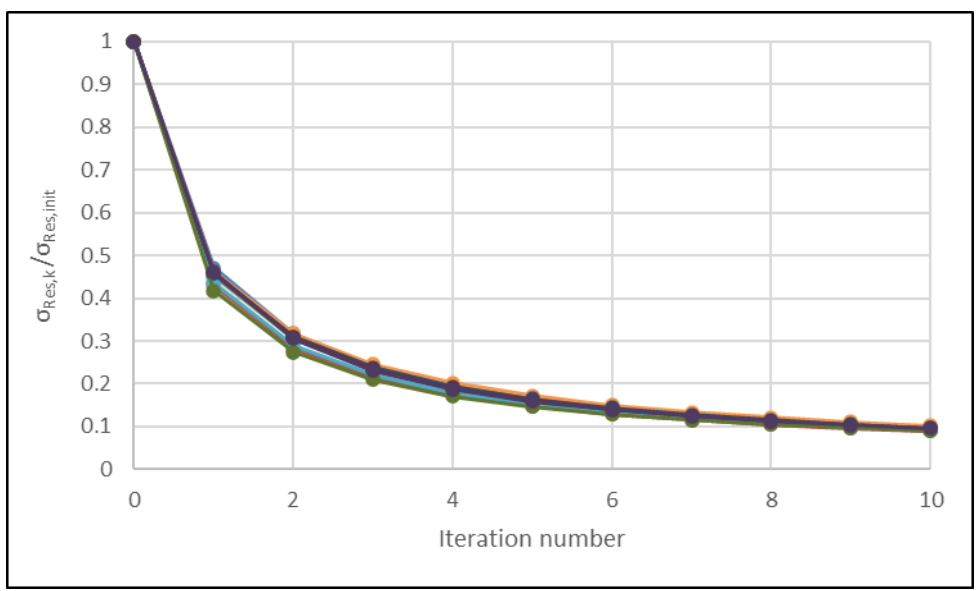

**Figure 4**. Residuals vs. iteration number for N = 1000 for 10 values of $x_0$

Figures 3 and 4 follow the same general pattern observed in Figure 1. Another observation that emerges from Figures 3 and 4 is that the variation due to randomly selected starting points is large for N = 10 but is much smaller for N = 1000. The variation is expected to be even smaller for higher values of N. Figures 5 and 6 show the results for distance from solution, obtained for 10 different starting points for N = 10, and for N = 1000.

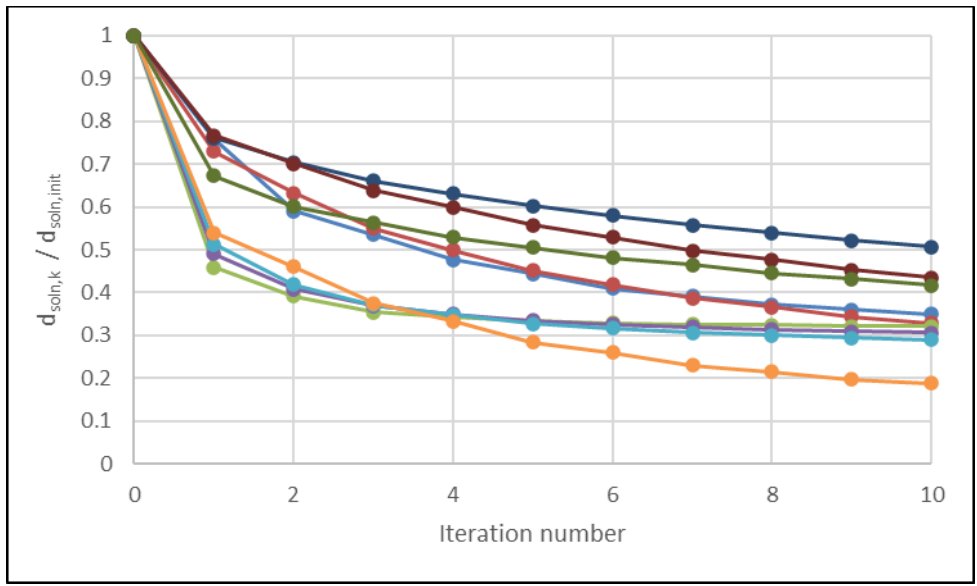

**Figure 5**. Distance from solution vs. iteration number for N = 10 for 10 values of $x_0$

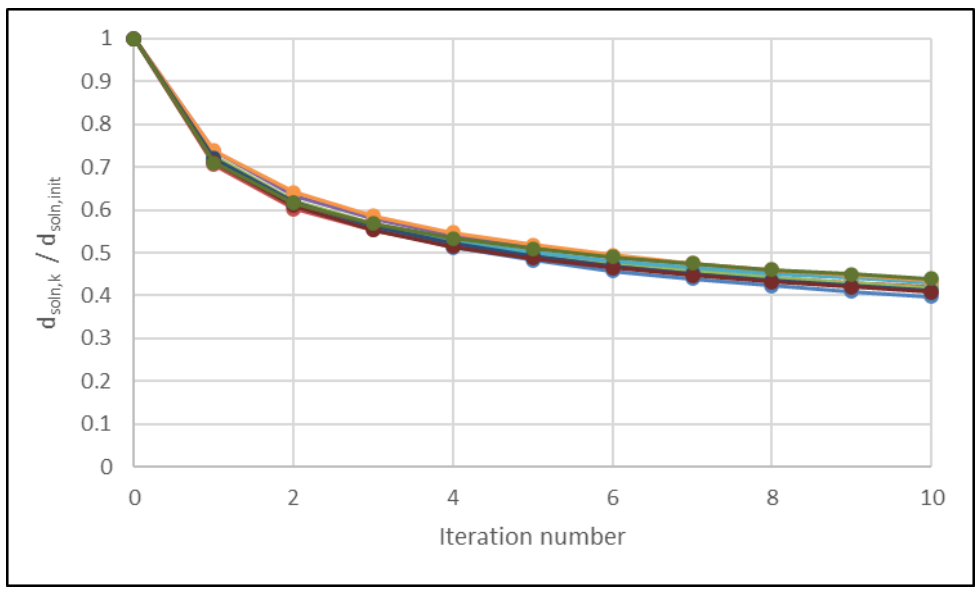

**Figure 6**. Distance from solution vs. iteration number for N = 1000 for 10 values of $x_0$

It is seen from these figures that both residuals and distance from solution reduce very fast in the first couple of steps, and then reduce further much more slowly in further steps. A natural question that arises is, whether we can modify or extend the SD procedure and get a faster convergence to the solution in later steps as well. If this can be done, the SD procedure can be used by itself to get very close to the solution in a reasonable number of steps.

In the following sections, two hypothetical (though impractical) possibilities are considered, which bring out the excellent potential performance of SD in approaching the solution.

## 4. Approach 1: Moving the iterate $x_k$ in a random direction before the next SD step

It has been shown that the progress made by successive SD iterates (in reducing the residuals or approaching the solution point) goes down rapidly with number of iterations. The lines shown in Figures 3 – 6 indicate that different initial points (selected at random) show similar behaviour: quick reduction of residuals, and fast approach to solution in the first iteration, i.e. by going from $x_0$ to $x_1$. However, further progress made by going from $x_1$ to $x_2$ is much lower. It is interesting to see the effect of moving from $x_1$ in a random direction, while preserving the distance from the solution, before taking the next step. The algorithm becomes:

______________________________________________

*SD, Algorithm 2: Moving X after each iteration*

Initial values:  $A$, $b$, and $x_0$ = the starting point
For k = 1 to n:  $r_k = b - A x_{k-1}$
  $\alpha_k = r_k^T r_k / ( r_k^T A r_k )$
  $x_k^{(1)} = x_{k-1} + \alpha_k r_k$
  $d_k$ = distance $[x_k^{(1)} - s]$
  $v$ = a random vector made with N(0,1) elements, then normalized
  $x_k = s + d_k v$

______________________________________________

It should be noted that the solution point **s** which has been used above, is unknown, and this algorithm has no practical use. However, the results of this algorithm serve the purpose of bringing out the excellent potential performance of SD. In practice, there may be other possible ways of moving the iterate $x_k$ in some (pseudo)random direction without changing the distance from the solution very much. This aspect will be considered in later publications.

Figure 1 shows the reduced residuals, $\sigma_{Res,k}/\sigma_{Res,init}$, for a specific problem for N = 10, up to 10 iterations. Figure 7 shows results obtained for the same problem up to 50 iterations, on a semi-log scale. Algorithm 1, which used SD alone, gives ($\sigma_{Res,k}/\sigma_{Res,init}$) of just 0.0054, while Algorithm 2, which moves the point $x_k$ after each iteration in a random direction, gives ($\sigma_{Res,k}/\sigma_{Res,init}$) of 2.6*(10$^{-8}$)! The difference is remarkable. Algorithm 2 thus gives a reduction in the residuals which is five orders of magnitude better than Algorithm 1. Figure 8 shows similar results for N = 1000.

Figures 7 and 8 show that Algorithm 2 reduces the residuals almost by a constant factor in each iteration. Thus, the residuals can be reduced to very low values in a reasonable number of iterations.

Similar conclusions can be drawn for distance from the solution. However, detailed results are not presented here, for the sake of brevity.

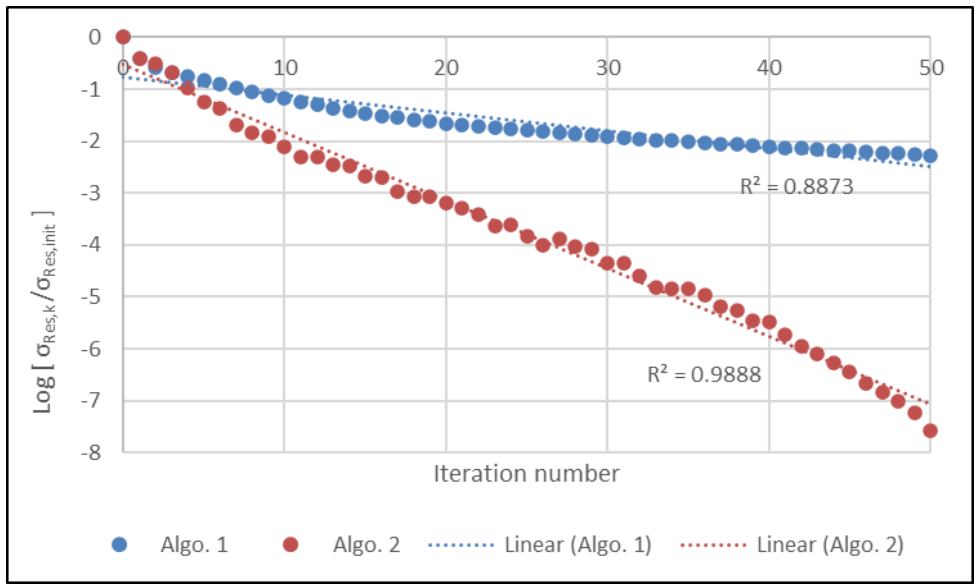

**Figure 7**. Residuals vs. iteration number for N = 10 for algorithms 1 and 2

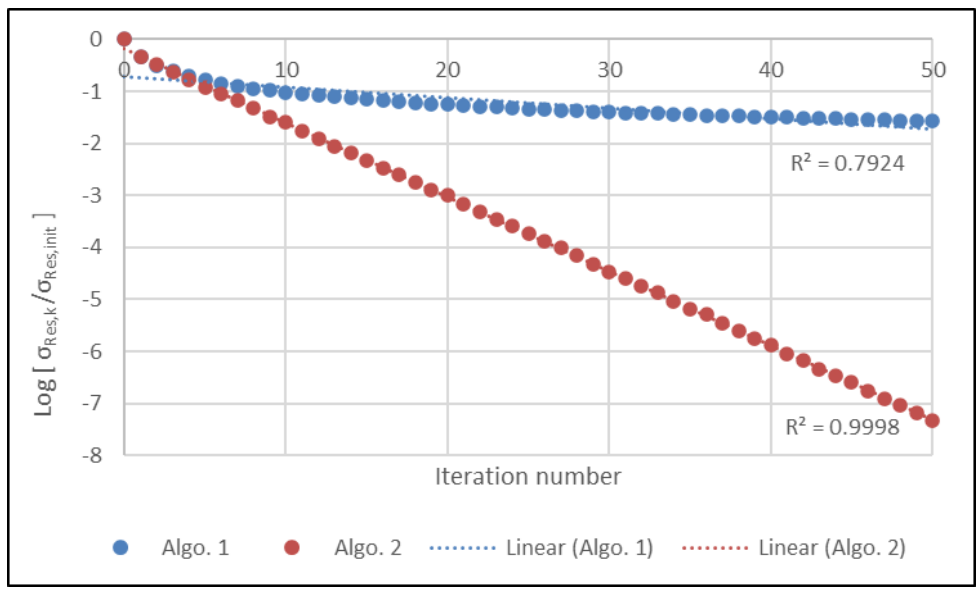

**Figure 8**. Residuals vs. iteration number for N = 1000 for algorithms 1 and 2

**5. Approach 2: Changing the matrix A in a random manner before the next SD step**

The system of equations, **Ax** = **b**, can be looked upon as N intersecting hyperplanes, all passing through the solution point **s**. Shifting the iterate $x_k$ to a new location before taking the next SD step, as described above, essentially changes the orientation of $x_k$ in relation to the intersecting hyperplanes. Another way of changing the orientation is to generate new instances of **A** and **b** at each iteration while keeping the same solution point **s**. (From a practical viewpoint, it is possible to transform **A** and **b** while maintaining the solution point, for example, by using row operations). The random generation of **A** at each iteration is an extreme case of such reorientations. Here we investigate the effect of this approach on the progress of iterations, to get some idea about the potential improvement in the performance of SD in approaching the solution. The algorithm becomes:

___________________________________________

*SD, Algorithm 3: Generating new **A** and **b**, while maintaining **s** after each iteration*

Initial values:         $A_0$, $b_0$ and $x_0$ (the starting point)
For k = 1 to n:         $r_k = b_{k-1} - A_{k-1}x_{k-1}$
                        $\alpha_k = r_k^T r_k / ( r_k^T A r_k )$
                        $x_k = x_{k-1} + \alpha_k r_k$
                        $A_k$ = a random matrix made with N(0,1) elements
                        $b_k = A_k s$

___________________________________________

It should be noted that the solution point **s** which has been used above, is unknown, and this algorithm has no practical use. However, the results of this algorithm serve the purpose of bringing out the excellent potential performance of SD. In practice, there may be other possible ways of transforming **A** and **b** at each iteration. This aspect will be considered in later publications.

Figure 1 shows the reduced residuals, $\sigma_{Res,k}/\sigma_{Res,init}$, for a specific problem for N = 10, up to 10 iterations. Figure 9 shows results obtained for the same problem up to 50 iterations, on a semi-log scale. Algorithm 1, which used just SD, gives ($\sigma_{Res,k}/\sigma_{Res,init}$) of just 0.0054, while Algorithm 3, which generates new **A** and **b** after each iteration in a random manner, gives ($\sigma_{Res,k}/\sigma_{Res,init}$) of $2.0*(10^{-10})$! The difference is remarkable. Algorithm 3 thus gives a reduction in the residuals which is seven orders of magnitude better than Algorithm 1. Figure 10 shows similar results for N = 1000.

Figure 9 and Figure 10 show that Algorithm 3 reduces the residuals almost by a constant factor in each iteration. Thus, the residuals can be reduced to very low values in a reasonable number of iterations. Similar conclusions can be drawn for distance from the solution. However, detailed results are not presented here, for the sake of brevity.

**6. Discussion of Results**

The results clearly show that SD leads to a fast approach to the solution (i.e. reduction in residuals, or distance from solution) in the first step itself, and slows down substantially in the following steps. The starting point used was one selected in a random direction from the solution. The extent of improvement obtained in the first step depends on the actual starting point and shows quite some variation for small sized problems (N = 10). However, for large sized problems (N = 1000) the variation seems to be very small.

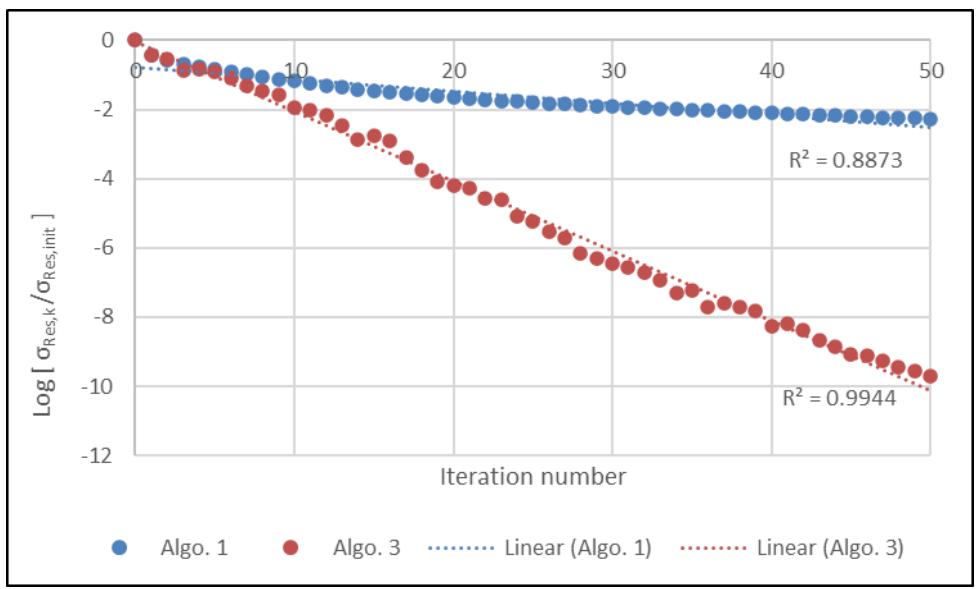

**Figure 9.** Residuals vs. iteration number for N = 10 for algorithms 1 and 3

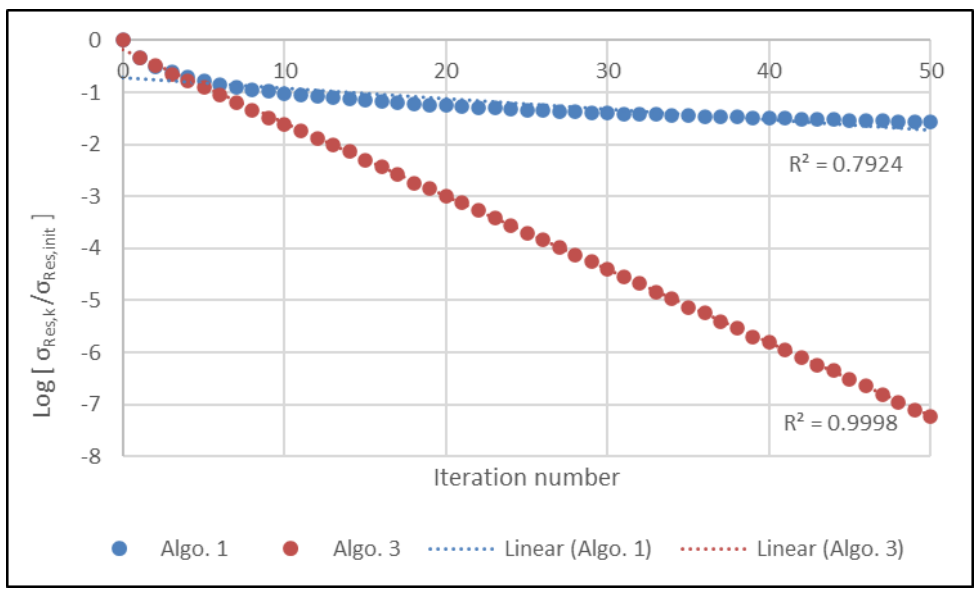

**Figure 10.** Residuals vs. iteration number for N = 1000 for algorithms 1 and 3

The slowing down of the SD progress with iterations can possibly be avoided, if the point obtained after an iteration is moved randomly before taking the next iteration. Two extreme approaches were investigated just to see how well such a strategy might work. The first approach involved moving the point in a random direction, while preserving the distance from the solution, before taking the next step. This essentially reorients the point in a random direction in relation to the **A** matrix. The details are given in Algorithm 2. The results in Figures 7 and 8 indicate that this approach can lead to a remarkable improvement in reduction of residuals, and the improvement can be by several orders of magnitude in 50 steps or so. Figure 7 shows that for small problems, such as N = 10, SD gives a smooth line, while Algorithm 2 gives some amount of scatter on the semi-log plot. The scatter almost disappears, as seen in Figure 8, for a large problem like N = 1000.

The second approach involved generating the **A** matrix and the **b** vector again in a random manner, while preserving the solution, before taking the next step. This essentially reorients and modifies the **A** matrix in relation to the current point. The details are given in Algorithm 3. The results in Figures 9 and 10 indicate that this approach can also lead to a remarkable improvement in reduction of residuals, and the improvement can be by several orders of magnitude in 50 steps or so. Figure 9 shows that for small problems, such as N = 10, SD gives a smooth line, while Algorithm 3 gives some amount of scatter on the semi-log plot. The scatter almost disappears, as seen in Figure 10, for a large problem like N = 1000.

The two approaches described above, and in Algorithms 2 and 3, use the solution point in moving the current point or in modifying **A** and **b**, before the next iteration step. This makes these approaches impractical to use, since the solution is not known beforehand. However, these results underline the tremendous improvement that similar approaches might lead to and make the SD method very effective. The development of such approaches, which do not use the solution point at all, will be the subject of later publications. The results presented here serve to set the stage for such developments.

## 7. Conclusions

The steepest descent method for solving simultaneous linear algebraic equations can be made substantially faster if we can adjust the position of an iterate or modify the **A** matrix and the **b** vector in a random manner. The exact method for doing such adjustments has not been presented here. However, the possibility of substantial improvement in speed of approaching the solution using SD (by several orders of magnitude) has been shown through computational experiments.